\documentclass[a4paper,12pt,reqno]{amsart}
\usepackage[a4paper,margin=1.9cm,top=2.5cm,bottom=2.5cm,centering,vcentering]{geometry}
\usepackage{amssymb,mathtools,cite,enumerate,color,eqnarray,hyperref}
\definecolor{blue}{RGB}{0, 0, 200}
\definecolor{pink}{RGB}{252, 0, 50}

\hypersetup{
           colorlinks=true,
           linkcolor=blue,
           urlcolor=red,
           citecolor=pink,
          }

\theoremstyle{plain}
\newtheorem{theorem}{Theorem}[section]

\newtheorem{lemma}[theorem]{Lemma}

\theoremstyle{definition}
\newtheorem{definition}[theorem]{Definition}
\newtheorem{Conjecture}[theorem]{Conjecture}

\numberwithin{equation}{section}

\numberwithin{equation}{section}

\begin{document}

\title[Congruences and density results for ped function]{Congruences and density results for partitions into distinct even parts}

\author[H. Nath]{Hemjyoti Nath}
\address[H. Nath]{Lokhra chariali, Guwahati 781040, Assam, India}
\email{hemjyotinath40@gmail.com}
\author[A. Sarma]{Abhishek Sarma}
\address[A. Sarma]{Department of mathematical sciences, Tezpur University, Napaam, Tezpur, Assam 784028, India}
\email{abhitezu002@gmail.com}
\keywords{Integer partitions, Ramanujan-type congruences, Radu's algorithm, arithmetic density.}

\subjclass[2020]{11P81, 11P83, 05A17, 11F11.}

\date{\today.}

\begin{abstract}
In this paper, we consider the set of partitions $ped(n)$ which counts the number of partitions of $n$ wherein the even parts are distinct (and the odd parts are unrestricted). Using an algorithm developed by Radu, we prove congruences modulo 192 which were conjectured by Nath \cite{21}. Further, we prove a few infinite families of congruences modulo 24 by using a result of Newman. Also, we prove that $ped(9n+7)$ is lacunary modulo $2^{k+2}\cdot 3$ and $3^{k+1}\cdot 4$ for all positive integers $k\geq0$. We further prove an infinite family of congruences for $ped(n)$ modulo arbitrary powers of 2 by employing a result of Ono and Taguchi on the nilpotency of Hecke operators.
\end{abstract}

\maketitle

\section{Introduction}
The partition of a positive integer $n$ is a non-increasing sequence of positive integers whose sum is equal to $n$. For example, 5+3+2 is a partition of 10.

If $p(n)$ denotes the number of partitions of a positive integer $n$, then with the convention that $p(0)=1$, the generating function of $p(n)$ (due to Euler) is given by
\begin{equation*}
    \sum_{n=0}^{\infty}p(n)q^n = \frac{1}{(q;q)_{\infty}},
\end{equation*}
where
\begin{equation*}
    (a;q)_{\infty} := \prod_{n=0}^{\infty}(1-aq^n),
\end{equation*}
where $a$ and $q$ are complex numbers with $|q|<1$. Throughout this paper, we set
\begin{equation*}
    f_k := (q^k;q^k)_{\infty}, \quad \text{for any integer} \quad k\geq1.
\end{equation*}

Ramanujan \cite{ram1}, \cite{ram2} and \cite{ram3} discovered three beautiful congruences for the partition function, namely
\begin{align*}
    p(5n+4)&\equiv 0 \pmod5, \\p(7n+5)&\equiv 0 \pmod7, \intertext{and} p(11n+7)&\equiv 0 \pmod{11}.
\end{align*}
The subject then got a lot of interest and many mathematicians over the years have found many more interesting results. Many Ramanujan-Type congruences for other classes of partitions such as $l$-regular partitions, etc (See for  example, \cite{9}) have also been found by many mathematicians. In this paper, the object of interest is the class of partitions where even parts are distinct and odd parts are unrestricted.

The number of partitions of $n$ wherein the even parts are distinct and the odd parts are unrestricted is denoted by $ped(n)$. With the convention that $ped(0)=1$, the generating function for $ped(n)$ is given by
\begin{equation}\label{e0}
    \sum_{n=0}^{\infty} ped(n)q^n = \frac{(-q^2;q^2)_{\infty}}{(q;q^2)_{\infty}} = \frac{(q^4;q^4)_{\infty}}{(q;q)_{\infty}}.
\end{equation} 
In \cite{1}, Andrews, Hirschhorn and Sellers investigated the $ped(n)$ function from an arithmetic point of view. In the same paper, they found several coungruences and exact generating functions for $ped(n)$. 
For example, they proved that
\begin{equation}\label{e2}
     \sum_{n=0}^{\infty}ped(9n+7)q^n = 12 \frac{f_2^4 f_3^6 f_4}{f_1^{11}}.
\end{equation}

Note that by $\eqref{e0}$, the number of partitions of $n$ wherein the even parts are distinct (and the odd parts are unrestricted) equals the number of partitions of $n$ with no part divisible by $4$, i.e., the $4$-regular partitions.,


In recent years, many arithmetic properties for the number of $4$-regular partitions have been studied. For related studies, one can see the following non-exhaustive list of  papers: Chen \cite{chendm}, Cui and Gu \cite{9}, Dai \cite{3}, Hirschhorn and Sellers \cite{hirschsel}, Merca \cite{mercajnt}, Xia \cite{7}, Du and Tang \cite{11.0} and Andrews and Bachraoui \cite{0.1}.

In a very recent work, Nath \cite{21} derived some congruences and an internal congruence for $ped(n)$ modulo $24$. In the same paper, he conjectured the following congruences for $ped(n)$.

\begin{Conjecture}\label{c1.0}
 For every $n\geq0$, we have
\begin{align*}
     ped(225n+43) & \equiv 0 \pmod{192},\\
     ped(225n+88) & \equiv 0 \pmod{192},\\
     ped(225n+133) & \equiv 0 \pmod{192},\\
     ped(225n+223) & \equiv 0 \pmod{192}.
\end{align*} 
\end{Conjecture}

The main aim of this paper is to prove the above conjecture, find some infinite families of congruences for $ped(n)$ modulo $24$ and prove the lacunarity of $\sum_{n=0}^{\infty} ped(9n+7)q^n$ modulo $2^{k+2}\cdot3$ and $3^{k+1}\cdot4$ for $k\geq0$. In the next four theorems, we state our main results.

\begin{theorem}\label{t1}
Conjecture $\ref{c1.0}$ is true.
\end{theorem}

For the next result, we require the concept of the Legendre symbol, which for a prime $p\geq3$ is defined as
\begin{align*}
\left(\dfrac{a}{p}\right)_L:=\begin{cases}\quad1,\quad \text{if $a$ is a quadratic residue modulo $p$ and $p\nmid a$,}\\\quad 0,\quad \text{if $p\mid a$,}\\~-1,\quad \text{if $a$ is a quadratic nonresidue modulo $p$.}
\end{cases}
\end{align*}
\begin{theorem}\label{t3}
Let $a(n)$ be defined by
    \begin{equation*}
        \sum_{n=0}^{\infty}a(n)q^n = f_1f_3^6.
    \end{equation*}
Let $p\geq5$ be a prime and let $\left(\frac{\star}{p}\right)_L$ denote the Legendre symbol. Define
\begin{equation*}
    \omega(p) := a\left( \frac{19}{24}(p^2-1) \right)+p^2\left( \frac{-2}{p} \right)_L\left( \frac{\frac{-19}{24}(p^2-1)}{p} \right)_L.
\end{equation*}
    \begin{enumerate}[(i)]
         \item If $\omega(p) \equiv 0 \pmod{2}$, then for $n,k\geq 0$ and $1\leq j \leq p-1$, we have
    \begin{equation}\label{t3.1}
            ped\left( 9p^{4k+4}n + 9p^{4k+3}j + \frac{57\cdot p^{4k+4}-1}{8} \right) \equiv 0 \pmod{24}.
        \end{equation}      
Furthermore, if $19 \nmid (24n+19)$, then for $n,k \geq 0$, we have
    \begin{equation}\label{t3.3.1}
            ped(9n+7) \equiv ped\left( 9\cdot 19^{4k+2}n + \frac{57\cdot 19^{4k+2}-1}{8} \right) \pmod{24}.
        \end{equation} 

        \item If $\omega(p) \equiv 1 \pmod{2}$, then for $n,k\geq 0$ and $1\leq j \leq p-1$, we have
      \begin{equation}\label{t2.2}
        ped\left( 9p^{6k+6}n + 9p^{6k+5}j + \frac{57\cdot p^{6k+6}-1}{8} \right) \equiv 0 \pmod{24}.
         \end{equation}
Furthermore, if $p \nmid (24n+19)$ and $p\neq 19$, then for $n,k\geq 0$ we have

    \begin{equation}\label{t3.3}
            ped\left( 9p^{6k+2}n + \frac{57\cdot p^{6k+2}-1}{8} \right) \equiv 0 \pmod{24}.
        \end{equation}
\end{enumerate}      
\end{theorem}

The distribution of the partition function modulo positive integers $M$ is also an interesting problem to explore. Given an integral power series $A(q) := \sum_{n=0}^{\infty}a(n)q^n $ and $0 \leq r \leq M$, we define the arithmetic density $\delta_r(A,M;X)$ as
\begin{equation*}
       \delta_r(A,M;X) = \frac{\#\{0 \leq n \leq X : a(n) \equiv r \pmod{M} \}}{X}.
  \end{equation*}
An integral power series A is called \textit{lacunary modulo $M$} if
\begin{equation*}
   \lim_{X \to \infty} \delta_0(A,M;X)=1,
\end{equation*}
that is, almost all the coefficients of $A$ are divisible by $M$.\\

The next theorem says that the partition function $ped(9n+7)$ is almost always divisible by $2^{k+2}\cdot3$ and $3^{k+1}\cdot4$ for $k\geq1$. To be specific, we prove the following result. 

\begin{theorem}\label{thm1}
	Let $G(q)=\sum_{n=0}^{\infty}ped(9n+7)q^n$. Then for every positive integer $k$, 
 \begin{equation}\label{e199}
     \lim_{X\to\infty} \delta_{0}(G,2^{k+2}\cdot 3;X)  = 1,
 \end{equation}
 \begin{equation}\label{e200}
     \lim_{X\to\infty} \delta_{0}(G,3^{k+1}\cdot 4;X)  = 1
 \end{equation}
	\end{theorem}

Also, we observe that the $eta$-quotient associated to $ped(9n+7)$ is a modular form whose level satisfies the conditions of Ono and Taguchi \cite[Theorem 1.3 (3)]{ot}. Thus, we use that result to prove the following congruences for $ped(9n+7)$.

\begin{theorem}\label{thm2}
Let $n$ be a nonnegative integer. Then there exists an integer $u\geq 0$
such that for every $v\geq 1$ and distinct primes $p_1,\ldots,p_{u+v}$ coprime to 6, we have
\begin{align*}
ped\left( \frac{3p_1\cdots p_{u+v}\cdot n-1}{8} \right)\equiv 0\pmod{2^{v}}
\end{align*}
whenever $n$ is coprime to $p_1,\ldots,p_{u+v}$.
\end{theorem}

The paper is organised as follows: In Section \ref{sec:pre}, we present some preliminaries required for our proofs. In Sections \ref{sec:t1}--\ref{proof:thm2}, we  present the proofs of Theorems \ref{t1}--\ref{thm2} respectively. Finally, in Section \ref{sec:conclude}, we end the paper with some concluding remarks.

\section{Preliminaries}\label{sec:pre}
The well known Ramanujan's general theta function $f(a,b)$ \cite[Equation 1.2.1]{6} is defined by
\begin{equation*}
    f(a,b)=\sum_{n=-\infty}^{\infty}a^{n(n+1)/2}b^{n(n-1)/2}, \quad |ab|<1.
\end{equation*}
Three special cases of $f(a,b)$ are the theta functions $\varphi(q)$, $\psi(q)$ and $f(-q)$, which are given by:
\begin{equation*}
    \varphi(q) := f(q,q) = \sum_{n=-\infty}^{\infty}q^{n^2} = (-q;q^2)_{\infty}^2(q^2;q^2)_{\infty} = \frac{f_2^5}{f_1^2f_4^2},
\end{equation*}
\begin{equation*}
    \psi(q) := f(q,q^3) = \sum_{n=0}^{\infty}q^{n(n+1)/2} = \frac{(q^2;q^2)_{\infty}}{(q;q^2)_{\infty}} = \frac{f_2^2}{f_1},
\end{equation*}
and
\begin{equation*}
    f(-q):=f(-q,-q^2)=\sum_{n=-\infty}^{\infty}(-1)^{n}q^{n(3n-1)/2} = (q;q)_{\infty}=f_1.
\end{equation*}
In terms of $f(a,b)$, Jacobi's triple product identity \cite[Equation 1.3.11]{6} is given by
\begin{equation*}
    f(a,b) = (-a;ab)_{\infty}(-b;ab)_{\infty}(ab;ab)_{\infty}.
\end{equation*}

The following result of Newman will play a crucial role in the proof of our second theorem, therefore we shall quote it as a lemma. Following the notations of Newman's paper, we shall let $p$, $q$ denote distinct primes, let $r, s \neq 0, r \not \equiv s \pmod{2}$. Set
\begin{equation}\label{e7}
    \phi(\tau) = \prod_{n=1}^{\infty}(1-x^n)^r(1-x^{nq})^s = \sum_{n=0}^{\infty}a(n)x^n,
\end{equation}
$\epsilon = \frac{1}{2}(r+s), t=(r+sq)/24, \Delta= t(p^2-1), \theta = (-1)^{\frac{1}{2}-\epsilon}2q^s$, then the result is as follows:

\begin{lemma}\label{l2}\textup{\cite[Theorem 3]{16}}
    With the notations defined as above, the coefficients $c(n)$ of $\phi(\tau)$ satisfy
    \begin{equation}\label{e8}
        a(np^2+\Delta)-\gamma(n)a(n) +p^{2\epsilon-2}a\left( \frac{n-\Delta}{p^2} \right) = 0,
    \end{equation}
where 
\begin{equation}\label{e9}
    \gamma(n) = p^{2\epsilon - 2}\alpha-\left( \frac{\theta}{p} \right)_L p^{\epsilon-3/2}\left( \frac{n-\Delta}{p} \right)_L,
\end{equation}
where $\alpha$ is a constant.
\end{lemma}


Now, we describe some useful background material on modular forms which will helpful for some of our proofs.

For a positive integer $N$, we assume that:
\begin{align*}
\textup{SL}_2(\mathbb{Z}) & :=\left\{\begin{bmatrix}
a  &  b \\
c  &  d      
\end{bmatrix}: a, b, c, d \in \mathbb{Z}, ad-bc=1
\right\},\\
\Gamma_{\infty} & :=\left\{
\begin{bmatrix}
1  &  n \\
0  &  1      
\end{bmatrix} \in \Gamma : n\in \mathbb{Z}  \right\},\\
\Gamma_{0}(N) &:=\left\{
\begin{bmatrix}
a  &  b \\
c  &  d      
\end{bmatrix} \in \Gamma : c\equiv~0\pmod N \right\},\\
\Gamma_{1}(N) &:=\left\{
\begin{bmatrix}
a  &  b \\
c  &  d      
\end{bmatrix} \in \Gamma : a\equiv d \equiv 1\pmod N \right\},\\
\Gamma(N) &:=\left\{
\begin{bmatrix}
a  &  b \\
c  &  d      
\end{bmatrix} \in \textup{SL}_2(\mathbb{Z}) : a \equiv d \equiv 1 \pmod{N}, \text{and} \hspace{2mm} b\equiv c \equiv 0 \pmod N \right\},\\
[\textup{SL}_2(\mathbb{Z}) : \Gamma_0(N)]&:=N\prod_{p\mid N} \left( 1+\dfrac{1}{p}\right),
\end{align*}
where $p$ is a prime number. A subgroup $\Gamma$ of $\textup{SL}_2(\mathbb{Z})$ is called a congruence subgroup if $\Gamma(N) \subseteq \Gamma$ for some $N$. The smallest $N$ such that $\Gamma(N) \subseteq \Gamma$ is called the level of $\Gamma$. For example, $\Gamma_{0}(N)$ and $\Gamma_{1}(N)$ are congruence subgroups of level $N$.\\

Let $\mathbb{H}:= \{ z \in \mathbb{C} : \Im(z) > 0 \}$ be the upper half of the complex plane. The group 

$$\textup{GL}_2^{+}(\mathbb{R}) = \left\{
\begin{bmatrix}
a  &  b \\
c  &  d      
\end{bmatrix} : a,b,c,d \in \mathbb{R} \hspace{2mm} \text{and} \hspace{2mm} ad-bc>0 \right\}$$ acts on $\mathbb{H}$ by $\begin{bmatrix}
a  &  b \\
c  &  d      
\end{bmatrix} z = \frac{az +b}{cz+d}$. We identify $\infty$ with $\frac{1}{0}$ and define $ \begin{bmatrix}
a  &  b \\
c  &  d      
\end{bmatrix} \frac{r}{s} = \frac{ar +bs}{cr+ds}$, where $\frac{r}{s} \in \mathbb{Q} \cup \{\infty\}$. This gives an action of $\textup{GL}_2^{+}(\mathbb{R})$ on the extended upper half-plane $\mathbb{H}^{\star} = \mathbb{H} \cup \mathbb{Q} \cup \{\infty\}$. Suppose that $\Gamma$ is a congruence subgroup of $\textup{SL}_2(\mathbb{Z})$. A cusp of $\Gamma$ is an equivalence class in $\mathbb{P}^{1}=\mathbb{Q} \cup \{ \infty\}$ under the action of $\Gamma$.\\

The group $\textup{GL}_2^{+}(\mathbb{R})$ also acts on the functions $f : \mathbb{H} \to \mathbb{C}$. In particular, suppose that $\gamma = \begin{bmatrix}
a  &  b \\
c  &  d      
\end{bmatrix} \in \textup{GL}_2^{+}(\mathbb{R}) $. If $f(z)$ is a meromorphic function on $\mathbb{H}$ and $\ell$ is an integer, then define the slash operator $\mid_{\ell}$ by 
$(f\mid_{\ell}\gamma )(z) := (det \gamma)^{\ell/2}(cz+d)^{-\ell}f(\gamma z)$.\\

\begin{definition}
Let $\Gamma$ be a congruence subgroup of level N. A holomorphic function $f : \mathbb{H} \to \mathbb{C}$ is called a modular form with integer weight $\ell$ on $\Gamma$ if the following hold:
\begin{enumerate}
    \item We have
    \begin{align*}
        f\left(\frac{az+b}{cz+d}\right)=(cz+d)^{\ell}f(z)
    \end{align*}
    for all $z\in\mathbb{H}$ and all $\begin{bmatrix}
        a&b\\
        c&d
    \end{bmatrix}\in\Gamma$.
    \item If $\gamma\in\textup{SL}_2(\mathbb{Z})$, then $(f\mid_{\ell}\gamma )(z)$ has a Fourier expansion of the form
    \begin{align*}
        (f\mid_{\ell}\gamma )(z)=\sum_{n=0}^{\infty}a_{\gamma}(n)q^n_{N},
    \end{align*}
    where $q^n_{N}=e^{\frac{2\pi iz}{N}}$.
\end{enumerate}
\end{definition}

For a positive integer $\ell$, let $M_\ell(\Gamma_1(N)))$ denote the complex vector space of modular forms of weight $\ell$ with respect to $\Gamma_1(N)$. 
\begin{definition}\cite[Definition 1.15]{ono2004}
If $\chi$ is a Dirichlet character modulo $N$, then a modular form $f\in M_\ell(\Gamma_1(N))$ has Nebentypus character $\chi$ if 
$f\left( \frac{az+b}{cz+d}\right)=\chi(d)(cz+d)^{\ell}f(z)$ for all $z\in \mathbb{H}$ and all $\begin{bmatrix}
a  &  b \\
c  &  d      
\end{bmatrix}\in \Gamma_0(N)$. The space of such modular forms is denoted by $M_\ell(\Gamma_0(N),\chi)$.
\end{definition}

To prove Theorem \ref{t1}, we  will need an algorithm developed by Radu. We discuss this next.
For integers $x$, let $[x]_m$ denote the residue class of $x$ in $\mathbb{Z}/m\mathbb{Z}$, $\mathbb{Z}^*_m$ denote the set of all invertible elements in $\mathbb{Z}_m$ and $\mathbb{S}_m$ denote the set of all squares in $\mathbb{Z}^*_m$. For integers $M\ge1$, let $R(M)$ be the set of all the integer sequences \[
(r_\delta):=\left(r_{\delta_1},r_{\delta_2},r_{\delta_3},\ldots,r_{\delta_k}\right)
\]
indexed by all the positive divisors $\delta$ of $M$, where $1=\delta_1<\delta_2<\cdots<\delta_k=M$. For integers $m\ge1$, $(r_\delta)\in R(M)$, and $t\in\{0,1,2,\ldots,m-1\}$, we define the set $P(t)$ as
\begin{align}
\label{Pt} P(t):=&\bigg\{t^\prime \in \{0,1,2,\ldots,m-1\} : t^{\prime}\equiv ts+\dfrac{s-1}{24}\sum_{\delta\mid M}\delta r_\delta \pmod{m} \notag\\
&\quad \text{~for some~} [s]_{24m}\in \mathbb{S}_{24m}\bigg\}.
\end{align}

For integers $m\ge1$, $M\ge1$, $N\ge1$,  $t\in \{0,1,2,\ldots,m-1\}$, $k:=\gcd\left(m^2-1,24\right)$, and $(r_{\delta})\in R(M)$, define $\Delta^{*}$ to be the set of all tuples $(m, M, N, t, (r_{\delta}))$ such that all of the following conditions are satisfied
	
	\begin{enumerate}
		\item[1.] Prime divisors of $m$ are also prime divisors of $N$;
		\item[2.] If $\delta\mid M$, then $\delta\mid mN$ for all $\delta\geq1$ with $r_{\delta} \neq 0$;
		\item[3.] $\displaystyle{24\mid kN\sum_{\delta\mid M}\dfrac{r_{\delta} mN}{\delta}}$;
		\item[4.] $\displaystyle{8\mid kN\sum_{\delta\mid M}r_{\delta}}$;
		\item[5.]  $\dfrac{24m}{\left(-24kt-k{\displaystyle{\sum_{\delta\mid M}}{\delta r_{\delta}}},24m\right)} \Bigm| N$;
  \item[6.] If $2|m$ then either $4|kN$ and $8|\delta N$ or $2|s$ and $8|(1-j)N$, where $\prod_{\delta|M}\delta^{|r_\delta|}=2^s\cdot j$.
	\end{enumerate}
Further, for integers $N\ge1$, $\gamma:=
\begin{pmatrix}
	a  &  b \\
	c  &  d
\end{pmatrix} \in \Gamma$, $(r_\delta)\in R(M)$, and $(r_\delta^\prime)\in R(N)$, we also define
	\begin{align*}
	p(\gamma)&:=\min_{\lambda\in\{0, 1, \ldots, m-1\}}\dfrac{1}{24}\sum_{\delta\mid M}r_{\delta}\dfrac{\gcd(\delta (a+ k\lambda c), mc)^2}{\delta m},\\
	p^\prime(\gamma)&:=\dfrac{1}{24}\sum_{\delta\mid N}r_{\delta}^\prime\dfrac{\gcd(\delta, c)^2}{\delta}.
	\end{align*}
 
The following lemma will support Lemma \ref{Lemma Radu 1} in the proof of Theorem \ref{t1}.
\begin{lemma}\label{Wang}\textup{\cite[Lemma 4.3]{Wang}} Let $N$ or $\frac{N}{2}$ be a square-free integer, then we have
		\begin{align*}
		\bigcup_{\delta\mid N}\Gamma_0(N)\begin{pmatrix}
		1  &  0 \\
		\delta  &  1
		\end{pmatrix}\Gamma_ {\infty}=\Gamma.
		\end{align*}
	\end{lemma}
 
We conclude this section with a result which will help us complete the proof of Theorem \ref{t1}.
\begin{lemma}\label{Lemma Radu 1}
\textup{\cite[Lemma 4.5]{Radu2}} Suppose that $(m, M, N, t, (r_{\delta}))\in\Delta^{*}$, $(r'_{\delta}):=(r'_{\delta})_{\delta \mid N}\in R(N)$, $\{\gamma_1,\gamma_2, \ldots, \gamma_n\}\subseteq \Gamma$ is a complete set of representatives of the double cosets of $\Gamma_{0}(N) \backslash \textup{SL}_2(\mathbb{Z})/ \Gamma_\infty$, and  $\displaystyle{t_{\min}:=\min_{t^\prime \in P(t)} t^\prime}$,
\begin{align}
	\label{Nu} \nu:= \dfrac{1}{24}\left( \left( \sum_{\delta\mid M}r_{\delta}+\sum_{\delta\mid N}r_{\delta}^\prime\right)[\textup{SL}_2(\mathbb{Z}):\Gamma_{0}(N)] -\sum_{\delta\mid N} \delta r_{\delta}^\prime-\frac{1}{m}\sum_{\delta|M}\delta r_{\delta}\right)
	- \frac{ t_{min}}{m},
 	\end{align}
$p(\gamma_j)+p^\prime(\gamma_j) \geq 0$ for all $1 \leq j \leq n$, and $\displaystyle{\sum_{n=0}^{\infty}A(n)q^n:=\prod_{\delta\mid M}f_\delta^{r_\delta}.}$ If for some integers $u\ge1$, all $t^\prime \in P(t)$, and $0\leq n\leq \lfloor\nu\rfloor$, $A(mn+t^\prime)\equiv0\pmod u$  is true,  then for integers $n\geq0$ and all $t^\prime\in P(t)$, we have $A(mn+t^\prime)\equiv0\pmod u$.
\end{lemma}

Recall that the Dedekind's eta-function $\eta(z)$ is defined by
\begin{align*}
	\eta(z):=q^{1/24}(q;q)_{\infty}=q^{1/24}\prod_{n=1}^{\infty}(1-q^n),
\end{align*}
where $q:=e^{2\pi iz}$ and $z\in \mathbb{H}$.

A function $f(z)$ is called an $eta$-quotient if it is expressible as a finite product of
the form $$f(z)=\prod_{\delta\mid N}\eta(\delta z)^{r_\delta},$$ where $N$ is a positive integer and each $r_{\delta}$ is an integer. The following two theorems allow one to determine whether a given eta-quotient is a modular form.
	\begin{theorem}\cite[Theorem 1.64]{ono2004}\label{thm_ono1} Suppose that $f(z)=\displaystyle\prod_{\delta\mid N}\eta(\delta z)^{r_\delta}$ 
		is an eta-quotient such that 
		\begin{eqnarray*}	
			\ell&=&\displaystyle\frac{1}{2}\sum_{\delta\mid N}r_{\delta}\in \mathbb{Z},\\
			\sum_{\delta\mid N} \delta r_{\delta}&\equiv& 0 \pmod{24} ~~\mbox{and}\\
			\sum_{\delta\mid N} \frac{N}{\delta}r_{\delta}&\equiv& 0 \pmod{24}.
		\end{eqnarray*}
		Then 
		$$
		f\left( \frac{az+b}{cz+d}\right)=\chi(d)(cz+d)^{\ell}f(z)
		$$
		for every  $\begin{bmatrix}
			a  &  b \\
			c  &  d      
		\end{bmatrix} \in \Gamma_0(N)$. Here 
		\begin{align}
		    \chi(d):=\left(\frac{(-1)^{\ell} \prod_{\delta\mid N}\delta^{r_{\delta}}}{d}\right)_L.\label{chi}
		\end{align}
	\end{theorem}

\noindent If the eta-quotient $f(z)$ satisfies the conditions of Theorem \ref{thm_ono1} and  holomorphic at all of the cusps of $\Gamma_0(N)$, then $f\in M_{\ell}(\Gamma_0(N), \chi)$. To determine  the orders of an eta-quotient at each cusp is the following.
	\begin{theorem}\cite[Theorem 1.65]{ono2004}\label{thm_ono1.1}
		Let $c, d,$ and $N$ be positive integers with $d\mid N$ and $\gcd(c, d)=1$. If $f(z)$ is an eta-quotient satisfying the conditions of Theorem~\ref{thm_ono1} for $N$, then the order of vanishing of $f(z)$ at the cusp $\frac{c}{d}$ 
		is $$\frac{N}{24}\sum_{\delta\mid N}\frac{\gcd(d,\delta)^2r_{\delta}}{\gcd(d,\frac{N}{d})d\delta}.$$
	\end{theorem}

\section{Proof of Theorem \ref{t1}}\label{sec:t1}
From \eqref{e2}, we see that to  prove Theorem \ref{t1}, it  is enough to prove
\begin{align}
        a(25n+4) & \equiv 0 \pmod{16},\label{cong1}\\
     a(25n+9) & \equiv 0 \pmod{16},\label{cong2}\\
     a(25n+14) & \equiv 0 \pmod{16},\label{cong3}\\
     a(25n+24) & \equiv 0 \pmod{16},\label{cong4}
\end{align}
where $\displaystyle\sum_{n=0}^{\infty}a(n)q^n=\frac{f_2^4 f_3^6 f_4}{f_1^{11}}\equiv \frac{f_1^5 f_3^6 f_4}{f_2^{4}} \pmod{16}$.

First, we prove \eqref{cong1} and \eqref{cong2}. We choose $(m,M,N,t,(r_\delta))=(25, 12, 60, 4, (5, -4, 6, 1, 0, 0))$. It is easy to check that this choice satisfies the required conditions of $\Delta^{\ast}$. By \eqref{Pt} we see that $P(t)=\{4, 9\}$. Using Lemma \ref{Wang}, we  see that $\left\{\left(\begin{pmatrix}
    1 &0\\
    \delta &1
\end{pmatrix}\right): \delta | N=60\right\}$ is a complete set of representatives of the double cosets in $\Gamma_{0}(N) \backslash \Gamma/ \Gamma_\infty$. Also, for the choice of $(r_\delta^\prime)=(0,0,0,0,0,0,0,0,0,0,0,0)$ we see that
\[p\left(\begin{pmatrix}
1 & 0\\
\delta & 1
\end{pmatrix}\right)+p^{\prime}\left(\begin{pmatrix}
1 & 0\\
\delta & 1
\end{pmatrix}\right)\ge 0\quad \text{for all $\delta\mid N=60$}.
\]
Using Lemma \ref{Lemma Radu 1}, we  find that
\begin{align*}
    \lfloor{\nu\rfloor}=53.
\end{align*} 
By choosing $u=16$ in Lemma \ref{Lemma Radu 1}, we need to check the congruences for all $n\leq \lfloor \nu \rfloor$, and then with the aid Lemma \ref{Lemma Radu 1} we can conclude our result. We verify this using Mathematica and hence the result follows.

For the proof of \eqref{cong3} and \eqref{cong4}, we choose $(m,M,N,t,(r_\delta))=(25, 12, 60, 14, (5, -4, 6, 1, 0, 0))$ and $(r_\delta^\prime)=(0,0,0,0,0,0,0,0,0,0,0,0)$. For $t=14$, we get $P(t)=\{14, 24\}$. The rest of the proof is exactly similar to the proof of the other two congruences, so we leave it to the reader.

Therefore, combining \eqref{cong1}, \eqref{cong2}, \eqref{cong3} and \eqref{cong4} with \eqref{e2}, we obtain Theorem \ref{t1}.

\section{Proof of theorem \ref{t3}}\label{sec:t3}
Again, from \eqref{e2}, we have that
    \begin{equation}\label{e10}
        \sum_{n=0}^{\infty}ped(9n+7)q^n = 12 \frac{f_2^4 f_3^6 f_4}{f_1^{11}} \equiv 12 f_1f_3^6 \equiv 12  \sum_{n=0}^{\infty}a(n)q^n \pmod{24},
    \end{equation}
    where $f_1f_3^6 = \displaystyle\sum_{n=0}^{\infty}a(n)q^n$.
    
    Putting $r=1, q=3$ and $s=6$ in $\eqref{e7}$, we have by Lemma $\ref{l2}$, for any $n\geq 0$
    \begin{equation}\label{e11}
        a\left( p^2n +\frac{19}{24}(p^2-1) \right) = \gamma(n)a(n) - p^5a\left( \frac{1}{p^2}\left( n-\frac{19}{24}(p^2-1) \right) \right)=0,
    \end{equation}
    where
    \begin{equation}\label{e12}
        \gamma(n)=p^5\alpha-p^2\left(\frac{-2}{p}\right)_L\left(\frac{n-\frac{19}{24}(p^2-1)}{p}\right)_L
    \end{equation}
    and $\alpha$ is a constant integer. Setting $n=0$ in $\eqref{e11}$ and using the fact that $a(0) = 1$ and $a\left(\frac{\frac{-19}{24}(p^2-1)}{p^2}\right)=0$, we obtain
    \begin{equation}\label{e13}
        a\left( \frac{19}{24}(p^2-1) \right) = \gamma(0).
    \end{equation}
Setting $n=0$ in $\eqref{e12}$ and using $\eqref{e13}$, we obtain
\begin{equation}\label{e14}
    p^5\alpha = a\left( \frac{19}{24}(p^2-1) \right)+p^2\left(
    \frac{-2}{p}\right)_L \left(\frac{\frac{-19}{24}(p^2-1)}{p}\right)_L := \omega(p).
    \end{equation}
Now rewriting $\eqref{e11}$, by referring $\eqref{e12}$ and $\eqref{e14}$, we obtain
\begin{align}
        a\left( p^2n +\frac{19}{24}(p^2-1) \right) &= \left( \omega(p)-p^2\left(\frac{-2}{p}\right)_L\left( \frac{n-\frac{19}{24}(p^2-1)}{p} \right)_L \right)a(n) \nonumber\\
        &\quad- p^5a\left( \frac{1}{p^2}\left( n-\frac{19}{24}(p^2-1) \right) \right).\label{e15}
    \end{align}
Now, replacing $n$ by $pn + \frac{19}{24}(p^2-1)$ in $\eqref{e15}$, we obtain
\begin{equation}\label{e16}
    a\left( p^3n+\frac{19}{24}(p^4-1) \right) = \omega(p)a\left( pn + \frac{19}{24}(p^2-1) \right) - p^5 a(n/p).
\end{equation}

From equation $\eqref{e11}$, we can see that
\begin{equation}\label{e28}
        a\left( p^2n +\frac{19}{24}(p^2-1) \right) - \gamma(n)a(n) + a\left( \frac{1}{p^2}\left( n-\frac{19}{24}(p^2-1) \right) \right) \equiv 0 \pmod{2},
    \end{equation}
    where
    \begin{equation}\label{e29}
        \gamma(n) \equiv \omega(p) + \left(\frac{n-\frac{19}{24}(p^2-1)}{p}\right)_L  \pmod2.
    \end{equation}
Setting $n=0$ in $\eqref{e28}$ and using the fact that $a(0)=1$ and $a\left( \frac{\frac{-19}{24}(p^2-1)}{p^2} \right) = 0$, we arrive at
\begin{equation}\label{e30}
    a\left( \frac{19}{24}(p^2-1) \right) \equiv \gamma(0) \pmod{2}.
\end{equation}
Setting $n=0$ in $\eqref{e29}$ yields
\begin{equation}\label{e31}
    \gamma(0) \equiv \omega(p) + 1 \pmod2, \qquad \text{for} \quad p \neq 19,
\end{equation}
and
\begin{equation}\label{e31.1}
    \gamma(0) \equiv \omega(p) \pmod2, \qquad \text{for} \quad p = 19.
\end{equation}
Combining $\eqref{e30}$ and $\eqref{e31}$ yields
\begin{equation}\label{e32}
    a\left( \frac{19}{24}(p^2-1) \right) + 1 \equiv \omega(p) \pmod{2},
\end{equation}
and combining $\eqref{e30}$ and $\eqref{e31.1}$ yields
\begin{equation}\label{e32:2}
    a\left( \frac{19}{24}(p^2-1) \right) \equiv \omega(p) \pmod{2}.
\end{equation}

\underline{\textbf{Case - 1 :} $\omega(p) \equiv 0 \pmod{2}$}

Since $\omega(p) \equiv 0 \pmod{2}$, from equation $\eqref{e16}$ we obtain
\begin{equation}\label{e17}
    a\left( p^3n+\frac{19}{24}(p^4-1) \right) \equiv p^5 a(n/p) \pmod{2}.
\end{equation}
Now, replacing $n$ by $pn$ in $\eqref{e17}$, we obtain
\begin{equation}\label{e18}
    a\left(p^4n+\frac{19}{24}(p^4-1)\right) \equiv p^5 a(n) \equiv a(n) \pmod{2}.
\end{equation}
Since $p^{4k}n+\frac{19}{24}(p^{4k}-1)=p^4\left(p^{4k-4}n+\frac{19}{24}(p^{{4k-4}}-1)\right)+\frac{19}{24}(p^4-1)$, using equation $\eqref{e17}$, we obtain that for every integer $k\geq 1$,
\begin{equation}\label{e19}
    a\left(p^{4k}n+\frac{19}{24}(p^{4k}-1)\right) \equiv a\left(p^{4k-4}n+\frac{19}{24}(p^{4k-4}-1)\right) \equiv a(n) \pmod{2}.
\end{equation}
Now if $p \nmid n$, then $\eqref{e17}$ yields
\begin{equation}\label{e20}
    a\left(p^3n +\frac{19}{24}(p^4-1)\right) \equiv 0 \pmod{2}.
\end{equation}
Replacing $n$ by $p^3+\frac{19}{24}(p^4-1)$ in $\eqref{e19}$ and using $\eqref{e20}$, we obtain
\begin{equation}\label{e21}
a\left( p^{4k+3}n + \frac{19}{24}(p^{4k+4}-1) \right) \equiv 0 \pmod{2}.    
\end{equation}
In particular, for $1 \leq j \leq p-1$, we have from $\eqref{e21}$, that
\begin{equation}\label{e40}
    a\left( p^{4k+4}n + p^{4k+3}j + \frac{19}{24}\left( p^{4k+4}-1 \right)\right) \equiv 0 \pmod{2}.
\end{equation}
Congruence $\eqref{t3.1}$ follows from \eqref{e10} and $\eqref{e40}$. \\

Now, assume $p=19$, we have from $\eqref{e29}$ and $\eqref{e32:2}$, that
\begin{equation}\label{e51}
    \gamma(n) \equiv a\left( \frac{19}{24}(p^2-1) \right) + \left( \frac{n-\frac{19}{24}(p^2-1)}{p} \right)_L \pmod{2}.
\end{equation}
By $\eqref{e28}$ and $\eqref{e51}$
\begin{align}
    a\left( p^2n + \frac{19}{24}(p^2-1) \right) &\equiv \left( a\left( \frac{19}{24}(p^2-1) \right) + \left( \frac{n - \frac{19}{24}(p^2-1)}{p} \right)_L \right) a(n) \nonumber \\
    &\quad + a\left( \frac{1}{p^2} \left( n - \frac{19}{24}(p^2-1) \right) \right) \pmod{2}. \label{e52}
\end{align}
Replacing $n$ by $pn + \frac{19}{24}(p^2-1)$ in $\eqref{e52}$ yields
\begin{equation}\label{e53}
    a\left( p^3n + \frac{19}{24}(p^4-1) \right) \equiv a\left( \frac{19}{24}(p^2-1) \right)a\left( pn + \frac{19}{24}(p^2-1) \right) + a(n/p) \pmod{2}.
\end{equation}
For $p=19$, we have $\omega(p) \equiv 0 \pmod{2}$ so from $\eqref{e32:2}$, we have $a\left( \frac{19}{24}(p^2-1) \right) \equiv 0 \pmod{2}$. Therefore, replacing $n$ by $pn$ in $\eqref{e53}$ we find that
\begin{equation}\label{e54}
    a\left( p^4n + \frac{19}{24}(p^4-1) \right) \equiv a(n) \pmod{2}.
\end{equation}
By $\eqref{e54}$ and iteration, we deduce that for $n,k \geq 0$,
\begin{equation}\label{e55}
    a\left( p^{4k}n + \frac{19}{24}(p^{4k}-1) \right) \equiv a(n) \pmod{2}.
\end{equation}
Moreover, we can rewrite $\eqref{e52}$ as 
\begin{equation}\label{e56}
    a\left( p^2n + \frac{19}{24}(p^2-1) \right) \equiv a(n) \left( \frac{n-\frac{19}{24}(p^2-1)}{p} \right) + a\left( \frac{1}{p^2}\left( n-\frac{19}{24}(p^2-1) \right) \right) \pmod{2}.
\end{equation}
If $p\nmid (24n+19)$, then $p \nmid \left( n -\frac{19}{24}(p^2-1) \right)$ and $\frac{n-\frac{19}{24}(p^2-1)}{p^2}$ is not an integer. Therefore
\begin{equation}\label{e57}
    \left( \frac{n-\frac{19}{24}(p^2-1)}{p} \right)_L \equiv 1 \pmod{2},
\end{equation}
and
\begin{equation}\label{e58}
    a\left( \frac{n-\frac{19}{24}(p^2-1)}{p^2} \right) = 0.
\end{equation}
On account of $\eqref{e56}$-$\eqref{e58}$, we deduce that
\begin{equation}\label{e59}
    a\left( p^2n + \frac{19}{24}(p^2-1) \right) \equiv a(n) \pmod{2}.
\end{equation}
Replacing $n$ by $p^2n + \frac{19}{24}(p^2-1)$ in $\eqref{e55}$ and employing $\eqref{e59}$, we see that for $k\geq 0$,
\begin{equation}\label{e60}
    a\left( p^{4k+2}n + \frac{19}{24}(p^{4k+2}-1) \right) \equiv a(n) \pmod{2}.
\end{equation}
Equations $\eqref{e10}$ and $\eqref{e60}$ readily yield $\eqref{t3.3.1}$.\\

\underline{\textbf{Case - 2 :} $\omega(p) \equiv 1 \pmod{2}$}

In order to prove $(ii)$, we replace $n$ by $p^2n+\frac{19}{24}p(p^2-1)$ in $\eqref{e16}$.
\begin{align}
    a\left( p^5n + \frac{19}{24}(p^6-1) \right) & = a\left( p^3\left( p^2n+\frac{19}{24}p(p^2-1) \right) + \frac{19}{24}(p^4-1) \right) \nonumber\\
    & \equiv \omega(p)a\left(p^3n+\frac{19}{24}(p^4-1)\right)-p^5a\left( pn + \frac{19}{24}(p^2-1) \right) \nonumber \\
    & \equiv \left[ \omega^2(p) - p^5\right]a\left( pn+\frac{19}{24}(p^2-1) \right) -p^5\omega(p)a(n/p). \label{e22}
\end{align}
Now, as $\omega(p) \equiv 1 \pmod{2}$ and $p\geq5$ is an odd prime, we have $\omega^2(p)-p^5=0 \pmod{2}$, and therefore $\eqref{e22}$ becomes
\begin{equation}\label{e23}
    a\left( p^5n + \frac{19}{24}(p^6-1) \right) \equiv a(n/p) \pmod{2}.
\end{equation}
Replacing $n$ by $pn$ in $\eqref{e23}$, we obtain
\begin{equation}\label{e24}
    a\left( p^6n +\frac{19}{24}(p^6-1) \right) \equiv a(n) \pmod{2}.
\end{equation}
Using equation $\eqref{e24}$ repeatedly, we see that for every integers $k\geq 1$,
\begin{equation}\label{e25}
    a\left( p^{6k+5}n + \frac{19}{24}(p^{6k}-1) \right) \equiv a(n) \pmod{2}.
\end{equation}
Observe that if $p\nmid n$, then $a(n/p)=0$. Thus $\eqref{e23}$ yields
\begin{equation}\label{e26}
    a\left( p^5n + \frac{19}{24}(p^6-1) \right) \equiv 0 \pmod{2}.
\end{equation}
Replacing $n$ by $p^5n+\frac{19}{24}(p^6-1)$ in $\eqref{e25}$ and using $\eqref{e26}$, we obtain
\begin{equation}\label{e27}
    a\left( p^{6k+5}n + \frac{19}{24}(p^{6k+6}-1) \right) \equiv 0 \pmod{2}.
\end{equation}
In particular, for $1 \leq j \leq p-1$, we have from $\eqref{e27}$, that
\begin{equation}\label{e27.1}
    a\left( p^{6k+6}n + p^{6k+5}j + \frac{19}{24}\left( p^{6k+6}-1 \right)\right) \equiv 0 \pmod{2}.
\end{equation}
Congruence $\eqref{t2.2}$ readily follows from $\eqref{e10}$ and $\eqref{e27.1}$. \\

Again, from $\eqref{e29}$  and $\eqref{e32}$, we have
\begin{equation}\label{e33}
    \gamma(n) \equiv a\left( \frac{19}{24}(p^2-1)\right) + 1 + \left( \frac{n-\frac{19}{24}(p^2-1)}{p} \right)_L \pmod{2}.
\end{equation}
By $\eqref{e28}$ and $\eqref{e33}$,
\begin{align}
    a\left( p^2n + \frac{19}{24}(p^2-1) \right) &\equiv \left( a\left( \frac{19}{24}(p^2-1) \right) + 1 + \left( \frac{n - \frac{19}{24}(p^2-1)}{p} \right)_L \right) a(n) \nonumber \\
    &\quad + a\left( \frac{1}{p^2} \left( n - \frac{19}{24}(p^2-1) \right) \right) \pmod{2} \label{e34}.
\end{align}

Replacing $n$ by $pn + \frac{19}{24}(p^2-1)$ in $\eqref{e34}$ yields 
\begin{equation}\label{e35}
    a\left( p^3n + \frac{19}{24}(p^4-1) \right) \equiv \left( a\left( \frac{19}{24}(p^2-1) \right) + 1 \right)a\left( pn + \frac{19}{24}(p^2-1) \right) + a(n/p).
\end{equation}
Since $\omega(p) \equiv 1 \pmod{2}$ so from $\eqref{e32}$, we have $a\left( \frac{19}{24}(p^2-1) \right) \equiv 0 \pmod{2}$. Therefore, replacing $n$ by $pn$ in $\eqref{e35}$, we find that 
\begin{equation}\label{e42}
    a\left( p^{4}n + \frac{19}{24}(p^4-1) \right) \equiv a\left( p^2n+\frac{19}{24}(p^2-1) \right) + a(n) \pmod{2}.
\end{equation}
Replacing $n$ by $p^2n + \frac{19}{24}(p^2-1)$ in $\eqref{e42}$ yields
\begin{equation}\label{e43}
    a\left( p^6n + \frac{19}{24}(p^6-1) \right) \equiv a\left( p^4n + \frac{19}{24}(p^4-1) \right) + a\left( p^2n+ \frac{19}{24}(p^2-1) \right) \pmod{2}.
\end{equation}
Combining $\eqref{e42}$ and $\eqref{e43}$ yields
\begin{equation}\label{e44}
    a\left( p^6n + \frac{19}{24}(p^6-1) \right) \equiv a(n) \pmod{2}.
\end{equation}
By $\eqref{e44}$ and iteration, we deduce that for $n,k \geq 0$,
\begin{equation}\label{e45}
    a\left( p^{6k}n + \frac{19}{24}(p^{6k}-1) \right) \equiv a(n) \pmod{2}.\\
\end{equation}

Moreover, we can rewrite $\eqref{e34}$ as 
\begin{equation}\label{e46}
    a\left( p^2n + \frac{19}{24}(p^2-1) \right) \equiv \left(1 + \left( \frac{n-\frac{19}{24}(p^2-1)}{p} \right)_L \right) a(n) + a\left( \frac{1}{p^2}\left( n-\frac{19}{24}(p^2-1) \right) \right) \pmod{2}.
\end{equation}
If $p\nmid (24n+19)$ and $p \neq 19$, then $p \nmid \left( n -\frac{19}{24}(p^2-1) \right)$ and $\frac{n-\frac{19}{24}(p^2-1)}{p^2}$ is not an integer. Therefore
\begin{equation}\label{e47}
    \left( \frac{n-\frac{19}{24}(p^2-1)}{p} \right)_L \equiv 1 \pmod{2},
\end{equation}
and
\begin{equation}\label{e48}
    a\left( \frac{n-\frac{19}{24}(p^2-1)}{p^2} \right) = 0.
\end{equation}
On account of $\eqref{e46}$-$\eqref{e48}$, we deduce that
\begin{equation}\label{e49}
    a\left( p^2n + \frac{19}{24}(p^2-1) \right) \equiv 0 \pmod{2}.
\end{equation}
Replacing $n$ by $p^2n + \frac{19}{24}(p^2-1)$ in $\eqref{e45}$ and employing $\eqref{e49}$, we see that for $k\geq 0$,
\begin{equation}\label{e50}
    a\left( p^{6k+2}n + \frac{19}{24}(p^{6k+2}-1) \right) \equiv 0 \pmod{2}.
\end{equation}
Congruence $\eqref{t3.3}$ follows from $\eqref{e10}$ and $\eqref{e50}$.\\

\section{Proof of Theorem \ref{thm1}}
 To prove Theorem \ref{thm1}, we need the following lemmas.
	\begin{lemma}\label{lem2}For any positive integers $k>2$,  we have 
		\begin{align}\label{e80}
			12 \frac{\eta(48z)^{4}\eta(72z)^6\eta(96z)\eta(24z)^{p^{k+1}-11}}{\eta(24pz)^{p^{k}}} \equiv \sum_{n=0}^{\infty}ped(9n+7)q^{24n+19} \pmod {p^{k+1}}.
		\end{align}
	\end{lemma}
	\begin{proof} 	
		Consider 
		\begin{align*}
			\mathcal{A}(z) 
			=\prod_{n=1}^{\infty}\frac{(1-q^{24n})^{p}}{(1-q^{24p n})}  =\frac{\eta(24z)^{p}}{\eta(24p z)}.
		\end{align*}
		By the binomial theorem, for any positive integers $r$, $k$, and prime $p$ we have
		\begin{align*}
			(q^{r};q^{r})_{\infty}^{p^k}\equiv (q^{pr};q^{pr})_{\infty}^{p^{k-1}} \pmod {p^{k}}.
		\end{align*}
		Therefore,
		\begin{align*}
			\mathcal{A}^{p^{k}}(z) = \frac{\eta(24z)^{p^{k+1}}}{\eta(24p z)^{p^k}} \equiv 1 \pmod {p^{k+1}}.
		\end{align*}
		Define $\mathcal{B}_{p, k}(z)$ by
		
		$$\mathcal{B}_{p, k}(z)=	\frac{\eta(48z)^4\eta(72 z)^6\eta(96z)}{\eta(24z)^{11}}~\mathcal{A}^{p^{k}}(z).$$
		
		\noindent Now, modulo $p^{k+1}$, we have 
		\begin{align}\label{new-110}
			\notag	\mathcal{B}_{p, k}(z)	&=\frac{\eta(48z)^4\eta(72 z)^6\eta(96z)}{\eta(24z)^{11}} \frac{\eta(24z)^{p^{k+1}}}{\eta(24p z)^{p^k}}\\
			\notag	&\equiv \frac{\eta(48z)^4\eta(72 z)^6\eta(96z)}{\eta(24z)^{11}}\\
			&= q^{19}\frac{f_{48}^4f_{72}^6f_{96}}{f_{24}^{11}}
		\end{align}
		Since
  \begin{equation}\label{e201}
    \mathcal{B}_{p, k}(z)=\frac{\eta(48z)^{4}\eta(72z)^6\eta(96z)\eta(24z)^{p^{k+1}-11}}{\eta(24p z)^{p^{k}}},
  \end{equation}
  combining \eqref{e2} and \eqref{new-110}, we obtain the required result.
	\end{proof}
	\begin{lemma}\label{lem1}
		For positive integers $k>2$, we have
		\begin{align*}
			\mathcal{B}_{2, k}(z) \in M_{2^{k-1}}\left(\Gamma_0(2304), \chi_1(\bullet)\right),
		\end{align*}
	where the Nebentypus character $$\chi_1(\bullet)=\left(\frac{(-1)^{2^{k-1}} \cdot 48^{4-2^k} \cdot 72^6 \cdot 96 \cdot 24^{2^{k+1}-11}}{\bullet}\right)_L.$$	
	\end{lemma}
	\begin{proof} We put $p=2$ in $\eqref{e201}$ to obtain

  \begin{equation*}
    \mathcal{B}_{2, k}(z)=\frac{\eta(72z)^6\eta(96z)\eta(24z)^{2^{k+1}-11}}{\eta(48 z)^{2^{k}-4}}.
  \end{equation*} 
 First we verify the first, second and third hypotheses of Theorem \ref{thm_ono1}. The weight of the $eta$-quotient $\mathcal{B}_{2, k}(z)$ is $2^{k-1}$.\\ Suppose the level of the $eta$-quotient $\mathcal{B}_{2, k}(z)$ is $288 u$, where $u$ is the smallest positive integer satisfying the following identity.
$$\frac{288 u}{12}+\frac{288 u}{96}+(2^{k+1}-11)\frac{288 u}{24}-(2^k-4)\frac{288 u}{48}\equiv0\pmod{24}.$$
Equivalently, we have
\begin{equation}
\label{ul}
    9u\left( 2^{k+1} - 9 \right)\equiv0\pmod{24}.
\end{equation}
Hence, from \eqref{ul}, we conclude that the level of the $eta$-quotient $\mathcal{B}_{2, k}(z)$ is $2304$ for $k>2$.

By Theorem \ref{thm_ono1.1}, the cusps of $\Gamma_{0}(2304)$ are given by  $\frac{c}{d}$ where $d~\mid~2304$ and $\gcd(c, d)~=~1$. Now  $\mathcal{B}_{2, k}(z)$ is holomorphic at a cusp $\frac{c}{d}$ if and only if
\begin{equation}\label{c1}	
    \mathcal{S}(2) := 8\frac{\gcd(d,72)^2}{\gcd(d,96)^2}+4(2^{k+1}-11)\frac{\gcd(d,24)^2}{\gcd(d,96)^2}-2(2^{k}-4)\frac{\gcd(d,48)^2}{\gcd(d,96)^2}+1\geq 0.
\end{equation}

It is not difficult to see that $\mathcal{S}(2)$ is positive. This is evident from the table below.

\begin{table}[h!]
\centering
\begin{tabular}{|c| c|}
\hline
$d \mid 2304$ & $\mathcal{S}(2)$  \\ \hline
$1, 2, 3, 4, 6, 8, 12, 24$ & $3\cdot 2^{k+1}-27$   \\
\hline
$9, 18, 36, 72$ & $3\cdot 2^{k+1}+37$   \\
\hline
$16, 48$ & $0$   \\
\hline
$144$ & $16$  \\
\hline
$32, 64, 96, 128, 192, 256, 384, 768$ & $3/4$   \\
\hline
$288, 576, 1152, 2304$ & $19/4$   \\ \hline
\end{tabular}

\end{table}

\begin{center}
    Table: Values of $\mathcal{S}(2)$ depending on $d$.
\end{center}
\noindent Therefore, from the above table, the orders of vanishing of $\mathcal{B}_{2, k}(z)$ at the at the cusp $\frac{c}{d}$ is nonnegative for $k>2$. So  $\mathcal{B}_{2, k}(z)$ is holomorphic at every cusp $\frac{c}{d}$. We have also verified the Nebentypus character  by Theorem \ref{thm_ono1}. Hence $\mathcal{B}_{2, k}(z)$ is a modular form of weight $2^{k-1}$ on $\Gamma_0(2304)$ with  Nebentypus character $\chi_1(\bullet)$. 
\end{proof}

	We state the following result of Serre, which is useful to prove  Theorem \ref{thm1}.
	\begin{theorem}\cite[Theorem~2.65]{ono2004}\label{serre}
		Let $k, m$  be positive integers. If  $f(z)\in M_{k}(\Gamma_0(N), \chi(\bullet))$ has the Fourier expansion $f(z)=\sum_{n=0}^{\infty}c(n)q^n\in \mathbb{Z}[[q]],$
		then there is a constant $\alpha>0$  such that
		$$
		\# \left\{n\leq X: c(n)\not\equiv 0 \pmod{m} \right\}= \mathcal{O}\left(\frac{X}{\log^{\alpha}{}X}\right).
		$$
	\end{theorem} 
	\begin{proof}[Proof of Theorem \ref{thm1}]Suppose $k>2$ is a positive integer. From Lemma~\ref{lem1}, we have 
		$$	\mathcal{B}_{2, k}(z)=\frac{\eta(48z)^{4}\eta(72z)^6\eta(96z)\eta(24z)^{2^{k+1}-11}}{\eta(48 z)^{2^{k}}} \in M_{2^{k-1}}\left(\Gamma_0(2304), \chi_1(\bullet)\right).$$  Also the Fourier coefficients of the $eta$-quotient $\mathcal{B}_{p, k}(z)$ are integers. So, by Theorem \ref{serre} and Lemma \ref{lem2}, we can find a constant $\alpha>0$ such that
		$$
		\# \left\{n\leq X:ped(9n+7)\not\equiv 0 \pmod{2^{k+2}\cdot 3} \right\}= \mathcal{O}\left(\frac{X}{\log^{\alpha}{}X}\right).
		$$
		Hence  $$\lim\limits_{X\to +\infty}\frac{	\# \left\{n\leq X:ped(9n+7)\equiv 0 \pmod{2^{k+2}\cdot 3} \right\}}{X}=1.$$
This allows us to prove the required divisibility by $2^{k+2}\cdot 3$ for all $k>2$ and since if a number is divisible by $2^{k+2}$, where $k>2$, then it is also divisible for lower powers of $k$, that is it is also divisible for $2^4\cdot 3$, $2^3\cdot 3$ and $2^2\cdot 3$ as well, hence this proves the required divisibility by $2^{k+2}\cdot 3$ for all $k\geq0$. This completes the proof of \eqref{e199}.\\

We next prove $\eqref{e200}$. We put $p=3$ in $\eqref{e201}$ to obtain
\begin{eqnarray*}    
    \mathcal{B}_{3, k}(z)=\frac{\eta(48z)^{4}\eta(96z)\eta(24z)^{3^{k+1}-11}}{\eta(72 z)^{3^{k}-6}}
\end{eqnarray*}
Now, $\mathcal{B}_{3, k}(z)$ is an $eta$-quotient with $N=2304$. As before, the cusps of $\Gamma_0(2304)$ are represented by fractions $c/d$, where $d \mid 2304$ and $\gcd(c,d)=1$. By Theorem $\ref{thm_ono1.1}$, $\mathcal{B}_{3, k}(z)$ is holomorphic at a cusp $c/d$ if and only if 

\begin{equation}\label{c2}	
    \mathcal{S}(3) := 8\frac{\gcd(d,48)^2}{\gcd(d,96)^2}+4(3^{k+1}-11)\frac{\gcd(d,24)^2}{\gcd(d,96)^2}-4(3^{k-1}-2)\frac{\gcd(d,72)^2}{\gcd(d,96)^2}+1\geq 0.
\end{equation}

We verify the positivity with the help of the following table.
\begin{table}[h!]
\centering
\begin{tabular}{|c| c|}
\hline
$d \mid 2304$  & $\mathcal{S}(3)$ \\ \hline
$1, 2, 3, 4, 6, 8, 12, 24$  & $32\cdot 3^{k-1}-27$   \\
\hline
$9, 18, 36, 72$  & $37$   \\
\hline
$16, 48$  & $8\cdot 3^{k-1}$   \\
\hline
$144$ & $16$  \\
\hline
$32, 64, 96, 128, 192, 256, 384, 768$ & $(8\cdot 3^{k-1}+3)/4$  \\
\hline
$288, 576, 1152, 2304$ & $19/4$  \\ \hline
\end{tabular}
\end{table}
\begin{center}
    Table: Values of $\mathcal{S}(3)$ depending on $d$.
\end{center}

From the above table we can see that $\mathcal{S}(3) \geq 0$ for all $d \mid 2304$. By Theorem $\ref{thm_ono1}$, $\mathcal{B}_{3, k}(z) \in M_{3^k}\left(\Gamma_0(2304), \chi_2(\bullet)\right)$, where $\chi_2$ is the associated Nebentypus character which can be easily calculated using \eqref{chi}. Using the same reasoning and $\eqref{lem2}$, we find that $ped(9n+7)$ is divisible by $3^{k+1}\cdot 4$ for almost all $n \geq 0$. This completes the proof of $\eqref{e200}$.

\end{proof}

\section{Proof of Theorem \ref{thm2}}\label{proof:thm2}
 We recall the following theorem of Ono and Taguchi \cite{ot} on the nilpotency of Hecke operators. 
\begin{theorem}\cite[Theorem 1.3 (3)]{ot}\label{thmot}
Let $n$ be a nonnegative integer and $k$ be a positive integer. Let $\chi$ be a quadratic Dirichlet character of conductor $9\cdot2^a$. Then there is an integer $c \geq 0$ such that for every $f(z) \in M_k(\Gamma_0(9\cdot2^a),\chi) \cap \mathbb{Z}[[q]]$ and every $t\geq 1$,
\begin{align*}
f(z)|T_{p_1}|T_{p_2}|\cdots|T_{p_{c+t}}\equiv 0 \pmod{2^t}
\end{align*}
whenever the primes $p_1, \ldots, p_{c+t}$ are coprime to 6.
\end{theorem}

Now, we apply the above theorem to the modular form $\mathcal{B}(z)$ to prove Theorem \ref{thm2}.

Taking $p=2$ in \eqref{e80}, we have
\begin{align}\label{e203}
 \mathcal{B}_{2, k}(z) \equiv 12^{-1} \cdot \sum_{n=0}^{\infty}ped(9n+7)q^{24n+19} \pmod {2^{k+1}},
\end{align}
which yields
\begin{align*}
 \mathcal{B}_{2, k}(z) := \sum_{n=0}^{\infty}B_{2,k}(n)q^n \equiv 12^{-1} \cdot \sum_{n=0}^{\infty}ped\left( \frac{3n-1}{8}\right)q^{n} \pmod {2^{k+1}}.
\end{align*}

Note that $\mathcal{B}_{2, k}(z)\in M_{2^{k-1}}\left(\Gamma_0(9\cdot 2^8),\chi_1\right)$, we find that there is an integer $u\geq 0$ such that for any $v\geq 1$,
\begin{align*}
\mathcal{B}_{2, k}(z)\mid T_{p_1}\mid T_{p_2}\mid\cdots\mid T_{p_{u+v}}\equiv 0\pmod{2^{v}}
\end{align*}
whenever $p_1,\ldots, p_{u+v}$ are coprime to 6. From the definition of Hecke operators, we have that if $p_1,\ldots, p_{u+v}$ are distinct primes and if $n$ is coprime to $p_1\cdots p_{u+v}$,
then
\begin{equation}\label{Fa22}
B_{2,k}(p_1\cdots p_{u+v} \cdot n) \equiv 0\pmod{2^{v}}.
\end{equation}

Combining $\eqref{e203}$ and $\eqref{Fa22}$, we complete the proof of the theorem.

\section{Concluding Remarks}\label{sec:conclude}
The authors believe that a lot more congruences should be true for the $ped$ function. Based on numerical evidences, we would like to end this section with the following conjecture:
\begin{Conjecture}
    Let $p\geq 5$ be a prime with $\left( \frac{-2}{p} \right)_{L}=-1$. Also, let $t$ be a positive integer with $(t,6) = 1$ and $p\mid t$. Then for all $n\geq 0$ and $1\leq j \leq p-1$, we have 
    \begin{equation*}
        ped\left( 9\cdot t^2n + \frac{9\cdot t^2j}{p} + \frac{57\cdot t^2-1}{8} \right) \equiv 0 \pmod{192}.
    \end{equation*}
\end{Conjecture}

\section*{acknowledgement}
The authors are grateful to Dr. Hirakjyoti Das. Discussions with Dr. Das were extremely helpful in proving Theorem \ref{t1}.  The second author was partially supported by an institutional fellowship for doctoral research from Tezpur University. He thanks the funding institution.\\

\textbf{Data availability} There is no data in this manuscript.

\section*{Declaration}

\textbf{Conflict of interest} The authors declare no competing interests.

\bibliographystyle{alpha}
\bibliography{ref}

\end{document}